\documentclass[sigconf,natbib=false]{acmart}

\usepackage[american]{babel}
\usepackage[T1]{fontenc}
\usepackage[utf8]{inputenc}
\usepackage[american]{isodate}

\usepackage[ruled,norelsize]{algorithm2e}

\usepackage{amsmath}
\usepackage{amsthm}
\usepackage{bm}

\usepackage{graphicx}
\usepackage{pgf, tikz}
\usepackage{pgfplots}
\pgfplotsset{compat=1.17}
\usepgfplotslibrary{statistics}
\usepackage{placeins}

\usepackage{subfig}
\usepackage{pgfplots}

\usepackage[datamodel=acmdatamodel,style=acmnumeric]{biblatex}
\usepackage[autostyle=true]{csquotes}
\usepackage{listingsutf8}
\usepackage{multirow}
\usepackage{orcidlink}
\usepackage{seqsplit}
\usepackage{siunitx}
\usepackage{xcolor}

\definecolor{sign}{HTML}{b02a2d}
\definecolor{direction}{HTML}{007900}
\definecolor{regime}{HTML}{8c399e}
\definecolor{characteristic}{HTML}{1f5dc2}
\definecolor{mantissa}{HTML}{636363}
\definecolor{error}{HTML}{BD002A}
\definecolor{cellbg}{HTML}{EDEDED}

\definecolor{p-sign}{HTML}{FF5454}
\definecolor{p-regime}{HTML}{CC9966}
\definecolor{p-regime-term}{HTML}{996633}
\definecolor{p-exponent}{HTML}{0080FF}
\definecolor{p-fraction}{HTML}{000000}

\definecolor{posit}{HTML}{b02a2d}
\definecolor{bfloat16}{HTML}{007900}
\definecolor{takum}{HTML}{1f5dc2}
\definecolor{float}{HTML}{636363}

\makeatletter
\def\lst@makecaption{%
  \def\@captype{table}%
  \@makecaption
}
\makeatother

\setcopyright{acmlicensed}
\copyrightyear{2025}
\acmYear{2025}
\acmDOI{XXXXXXX.XXXXXXX}

\acmConference[SC '25]{International Conference for High Performance Computing, Networking, Storage and Analysis}{November 16--21,
  2025}{St. Louis, MO}

\acmBooktitle{SC '25: Proceedings of the International Conference for High Performance Computing, Networking, Storage and Analysis, November 16--21, 2025, St. Louis, MO}
\acmISBN{978-1-4503-XXXX-X/2025/11}

\begin{document}

\title{Numerical Performance of the Implicitly Restarted Arnoldi Method in OFP8, Bfloat16, Posit, and Takum Arithmetics}

\author{Laslo Hunhold}
\email{hunhold@uni-koeln.de}
\orcid{0000-0001-8059-0298}
\affiliation{%
  \institution{University of Cologne}
  \city{Cologne}
  \country{Germany}
}

\author{James Quinlan}
\email{james.quinlan@maine.edu}
\orcid{0000-0001-8059-0298}
\affiliation{%
  \institution{University of Southern Maine}
  \city{Portland}
  \state{Maine}
  \country{USA}
}

\author{Stefan Wesner}
\email{wesner@uni-koeln.de}
\orcid{0000-0002-7270-7959}
\affiliation{%
  \institution{University of Cologne}
  \city{Cologne}
  \country{Germany}
}

\begin{abstract}
The computation of select eigenvalues and eigenvectors of large, sparse matrices is fundamental to a wide range of applications. Accordingly, evaluating the numerical performance of emerging alternatives to the IEEE 754 floating-point standard---such as OFP8 (E4M3 and E5M2), \texttt{bfloat16}, and the tapered-precision posit and takum formats---is of significant interest. Among the most widely used methods for this task is the implicitly restarted \textsc{Arnoldi} method, as implemented in ARPACK.
\par
This paper presents a comprehensive and untailored evaluation based on two real-world datasets: the SuiteSparse Matrix Collection, which includes matrices of varying sizes and condition numbers, and the Network Repository, a large collection of graphs from practical applications. The results demonstrate that the tapered-precision posit and takum formats provide improved numerical performance, with takum arithmetic avoiding several weaknesses observed in posits. While \texttt{bfloat16} performs consistently better than \texttt{float16}, the OFP8 types are generally unsuitable for general-purpose computations.
\end{abstract}




\maketitle
\section{Introduction}
The eigenvalue problem of finding an eigenpair $(\lambda,x)$ satisfying $Ax=\lambda x$ for a given real matrix $A$ is foundational in scientific computing, with a broad range of applications. In addition to its central role in numerical analysis for obtaining singular values and eigenvalues in general matrix computations, it is particularly significant in graph theory. Here, the eigenvalues of (typically sparse and large-scale) Laplacian or adjacency matrices encode crucial structural information, including graph connectivity \cite{Atay2020}, random walks \cite{Lovasz1993}, transportation problems \cite{munkres1957algorithms}, community structure \cite{Newman2006}, and the chromatic number \cite{Wilf1967}. Eigenvalues and eigenvectors of graph-associated matrices also reveal properties such as node centrality and network influence.
\par
A widely used algorithm for computing selected eigenvalues and eigenvectors of large sparse matrices is the implicitly restarted \textsc{Arnoldi} method. This algorithm belongs to the family of Krylov subspace methods and is implemented in the ARPACK software library \cite{arpack}.
\par
As the size of problems under investigation continues to grow and processor speed outpaces memory bandwidth, communication has become the primary bottleneck in high-performance computing---a phenomenon often referred to as the \enquote{memory wall.} To address this, there is increasing interest in using lower-precision arithmetic to reduce memory traffic. However, the IEEE 754 standard for floating-point arithmetic proves inadequate for low-precision scenarios due to its limited numerical accuracy and dynamic range at reduced bit widths.
\par
In response, several alternative number formats have been proposed. At 8 bits, the OFP8 formats E4M3 and E5M2 have emerged \cite{ofp8}. The E4M3 format allocates 4 bits for the exponent and 3 bits for the fraction, while E5M2 uses 5 exponent bits and 2 fraction bits—providing a wider dynamic range at the cost of reduced precision. In addition to these, \texttt{bfloat16}---a format proposed by Google Brain and now a widely adopted alternative to \texttt{float16} \cite{bfloat16}---has gained traction. These formats have recently been implemented in Intel’s AVX10.2 vector instruction set architecture \cite{intel-avx10.2}.
\par
While the aforementioned formats adhere to the IEEE 754 framework, more radical alternatives have been explored. One such approach is posit arithmetic, a tapered-precision format introduced by Gustafson et al. in 2017 \cite{posits-beating_floating-point-2017}. In posits, a variable-length exponent encoding allocates more fraction bits to values near one in magnitude and fewer to those farther away, mirroring the statistical distribution of values typically encountered in computation.
\par
Building on the posit framework, takum arithmetic was proposed in 2024 \cite{2024-takum}. Takums offer a significantly wider dynamic range and improved precision for values of very large or small magnitude, albeit with slightly reduced precision near one. They also aim to simplify hardware implementation. In this study, we focus on the floating-point variant known as linear takums, which we refer to simply as takums throughout the paper.
\par
This paper evaluates the numerical performance of emerging arithmetic formats through a series of experiments based on the implicitly restarted \textsc{Arnoldi} method, applied to sparse, symmetric matrices originating from both general computational problems (specifically from the SuiteSparse Matrix Collection \cite{davis2011university}) and graph Laplacians derived from large-scale networks in the Network Repository \cite{network_repository}. The algorithms used are not tailored to specific formats to preserve fairness and ensure general applicability. While this study does not include an evaluation of execution time or energy consumption---since the formats are mostly implemented in software---its findings remain highly relevant for practitioners. Numerical accuracy represents an upper bound on the capabilities of any arithmetic format, whereas performance and energy efficiency are variable targets subject to future optimization.
\par
To the best of the authors' knowledge, this is the first comprehensive study to systematically evaluate the use of OFP8, \texttt{bfloat16}, posits, and takums in the context of eigenvalue computation, particularly with the \textsc{Arnoldi} method. This work presents a unified evaluation across all these formats using a common experimental framework. In short, it constitutes the first broad assessment of emerging number formats for eigenvalue problems.
\par
The remainder of the paper is organized as follows: Section~\ref{sec:methods} details the experimental setup and benchmarking methodology. Section~\ref{sec:results} presents the main results, including analyses and visualizations. Section~\ref{sec:conclusions} concludes the paper with a summary of findings and final remarks.
\section{Experimental Methods}\label{sec:methods}
The experimental framework employed in this study builds upon the MuFoLAB software package \cite{mufolab}, which provides a systematic and reproducible environment for evaluating large-scale benchmarks. MuFoLAB offers a generic interface for conducting experiments, along with automated data preparation and analysis capabilities. A detailed framework description is available in \cite[Section~2]{mufolab-paper}.
\par
For the purposes of this study, the framework has been extended to support OFP8 numerical types, include preprocessing routines for graph datasets to produce sparse Laplacian matrices, and incorporate an experimental interface for eigenvalue computation via the \textsc{Arnoldi} method. The preparation of the 302 general, symmetric matrices sourced from the SuiteSparse Matrix Collection follows the procedure outlined in \cite{mufolab-paper}, with the only modification being the removal of all non-symmetric matrices and those containing more than 20,000 non-zero entries from the original set. The following subsection outlines the newly introduced components and the corresponding design decisions. Based on these components, the individual experiments are implemented in \cite[src/eigen\_*]{mufolab}.
\subsection{Graph Laplacian Matrices Preparation}
The graph datasets used in this study were obtained from the Network Repository \cite{network_repository}, following a three-stage process: scraping, cleanup with preprocessing and compression, and subsequent aggregation into four categories. A key challenge was that the Network Repository does not provide a single downloadable archive or a public API. Instead, users are expected to manually download individual datasets through the website interface—a method unsuitable for the scale of this study. To address this, the retrieval process was fully automated.
\par
The first stage, scraping, is implemented via a shell script \cite[src/generate\_graphs.jl]{mufolab}. This script downloads the HTML content of the Network Repository’s central listing page at \path{https://networkrepository.com/networks.php}, and uses \texttt{sed(1)} to extract all graph ZIP file download links, which consistently follow the pattern \path{https://nrvis.com/./download/data/}, followed by the respective graph identifier of the form \texttt{category/name}. This pattern enables automated extraction of all available download links. Internally, rather than storing a list of direct URLs, the system maintains a list of the 5,158 graph identifiers.
\par
To ensure data quality and manageability, a blacklist of 61 manually identified malformed entries is first applied to exclude problematic files. Next, a file size filter limits the dataset to graphs whose corresponding ZIP archives do not exceed \SI{500}{\kilo\byte}. To avoid two HTTP requests per graph---one to check the file size and another to download--the script employs \texttt{curl(1)} with the \texttt{--max-filesize} flag. If a file exceeds the limit, the download is aborted and the graph is discarded based on \texttt{curl}'s return code. After these two filtering steps, 3,302 graph archives remain. Each archive is unpacked and placed into the output directory \texttt{out/graphs/}, sorted into subdirectories corresponding to the original category of the graph.
\par
The Network Repository provides graphs in two formats: Matrix Market files (with a \texttt{.mtx} extension) and edge list files (with a \texttt{.edges} extension). Due to inconsistencies in the formatting of these files, a post-processing step is required to standardize them for further use. This includes applying both general parsing rules and file-specific adjustments. The final unpacked dataset occupies \SI{680}{\mega\byte} across 31 categories (see Table~\ref{tab:classes}) within the \texttt{out/graphs} directory. For reproducibility and integration with the \texttt{make(1)} build system, the script also generates a witness file, \texttt{src/generate\_graphs.output}, which lists all successfully processed graph files.
\par
The second stage---cleanup with preprocessing and compression---begins by reading the witness file generated in the previous step. Each listed graph is then loaded using a separate Julia program \cite[src/generate\_graph\_matrices.jl]{mufolab}. Matrix Market files are parsed with the \texttt{MatrixMarket.jl} package, while edge list files are read using \texttt{GraphIO.EdgeList}. As all matrices are intended to represent adjacency matrices, a necessary structural condition is squareness. However, several files violate this constraint due to extraneous zero blocks. These cases are identified by inspecting matrix dimensions; if a zero block can be removed to produce a square matrix, it is discarded. Conversely, if no such zero block is present, we instead append a zero block to ensure that the matrix is square.
\par
The next step is the computation of the symmetrically normalized Laplacian. The input matrix is first symmetrized via average symmetrization, defined by the transformation $A \mapsto \frac{1}{2}(A + A^T)$. This standard technique in graph spectral analysis ensures the resulting graph is undirected and is also required for our analysis, as discussed later. Given a symmetric adjacency matrix $A$, with vertex degrees defined as $\deg(i) := \sum_j A_{ij}$, the symmetrically normalized Laplacian is constructed as
\begin{equation}
L^\text{sym}_{i,j} := \begin{cases}
    1 & i = j \land \deg(i) > 0 \\
    -\frac{1}{\sqrt{\deg(i)\deg(j)}} & i \neq j \land A_{i,j} \neq 0 \\
    0 & \text{otherwise.}
\end{cases}
\end{equation}
By construction, $L^\text{sym}$ is symmetric and has unit diagonal entries for non-isolated vertices, effectively normalizing all non-zero degrees to one.
\par
Each Laplacian matrix is then wrapped in a \texttt{TestMatrix} structure that stores metadata and the matrix itself. These objects are collected into an array, sorted lexicographically by name, and serialized into a compressed \texttt{JLD2} file, \path{out/graph_test_matrices.jl}, with a final size of \SI{257}{\mega\byte}. This format ensures portability and efficient loading of the matrices within Julia.
\par
\begin{table}[tbp]
	\begin{center}
		\bgroup
		\def\arraystretch{1.3}
		\setlength{\tabcolsep}{0.3em}
        \begin{tabular}{|l|l|l|l|}
        \hline
        \textbf{class}                  & \textbf{class size}   & \textbf{graph category} & \textbf{category size} \\ \hline\hline
        \multirow{4}{*}{biological}     & \multirow{4}{*}{1219} & bio                     & 24                     \\ \cline{3-4} 
                                        &                       & eco                     & 6                      \\ \cline{3-4} 
                                        &                       & protein                 & 1178                   \\ \cline{3-4} 
                                        &                       & bn                      & 11                     \\ \hline\hline
        \multirow{6}{*}{infrastructure} & \multirow{6}{*}{29}   & inf                     & 4                      \\ \cline{3-4} 
                                        &                       & massive                 & 0                      \\ \cline{3-4} 
                                        &                       & power                   & 8                      \\ \cline{3-4} 
                                        &                       & road                    & 3                      \\ \cline{3-4} 
                                        &                       & tech                    & 5                      \\ \cline{3-4} 
                                        &                       & web                     & 9                      \\ \hline\hline
        \multirow{13}{*}{social}        & \multirow{13}{*}{23} & ca                      & 7                      \\ \cline{3-4} 
                                        &                       & cit                     & 1                      \\ \cline{3-4} 
                                        &                       & dynamic                 & 43                     \\ \cline{3-4} 
                                        &                       & econ                    & 12                     \\ \cline{3-4} 
                                        &                       & email                   & 6                      \\ \cline{3-4} 
                                        &                       & ia                      & 17                     \\ \cline{3-4} 
                                        &                       & proximity               & 6                      \\ \cline{3-4} 
                                        &                       & rec                     & 2                      \\ \cline{3-4} 
                                        &                       & retweet\_graphs         & 28                     \\ \cline{3-4} 
                                        &                       & rt                      & 31                     \\ \cline{3-4} 
                                        &                       & soc                     & 21                     \\ \cline{3-4} 
                                        &                       & socfb                   & 27                     \\ \cline{3-4} 
                                        &                       & tscc                    & 33                     \\ \hline\hline
        \multirow{8}{*}{miscellaneous}  & \multirow{8}{*}{1820} & dimacs                  & 62                     \\ \cline{3-4} 
                                        &                       & dimacs10                & 17                     \\ \cline{3-4} 
                                        &                       & graph500                & 0                      \\ \cline{3-4} 
                                        &                       & heter                   & 0                      \\ \cline{3-4} 
                                        &                       & labeled                 & 47                     \\ \cline{3-4} 
                                        &                       & misc                    & 1555                   \\ \cline{3-4} 
                                        &                       & rand                    & 139                    \\ \cline{3-4} 
                                        &                       & sc                      & 0                      \\ \hline
        \end{tabular}
		\egroup
	\end{center}
	\caption{Classification of graphs into four classes, composed from all 31 categories in the Network Repository. The counts reflect the 500\,kB size limit of the compressed graph archives, which accounts for some categories being empty.}
	\label{tab:classes}
\end{table}
The final stage---aggregation into four categories--- occurs during the loading phase \cite[src/TestMatrices.jl]{mufolab}. Here, the original 31 categories are aggregated into four broader classes: \emph{biological}, \emph{infrastructure}, \emph{social}, and \emph{miscellaneous} (see Table~\ref{tab:classes} for a detailed mapping). These four classes serve as the basis for the subsequent analysis presented in this work.
\subsection{Eigenvalue Experiment Interface}
With both the general matrices and graph Laplacians prepared, we now describe the execution of the eigenvalue experiments themselves. MuFoLAB provides a unified experimental interface \cite[Section~II.B]{mufolab}, which encapsulates all functionality for parallel computation and postprocessing of experiment results \cite[src/Experiments.jl]{mufolab}. This section assumes that each input matrix has already been converted to the numeric format under evaluation.
\par
At the core of the eigenvalue experiments lies the \textsc{Arnoldi} method with \textsc{Krylov}-\textsc{Schur} restarts---commonly known as the implicitly restarted \textsc{Arnoldi} method. This algorithm is implemented in the Julia package \texttt{ArnoldiMethod.jl}, which offers a type-generic interface suitable for a wide range of matrix representations.
\par
We employ the function \texttt{partialschur()} to compute a partial \textsc{Schur} decomposition of the input matrix. This function targets a user-specified number of eigenvalues, given by \texttt{eigenvalue\_count}, and accepts an ordering rule as input (e.g., to return either the smallest or largest eigenvalues). In all experiments, we configure the algorithm to return the 10 largest eigenvalues.
An additional parameter to \texttt{partialschur()} is the relative convergence tolerance. We set this tolerance to $10^{-2}$ for all 8-bit types, $10^{-4}$ for 16-bit types, $10^{-8}$ for 32-bit types, $10^{-12}$ for 64-bit types, and $10^{-20}$ when computing reference solutions in \texttt{float128}.
\par
For symmetric input matrices $M$, the resulting \textsc{Schur} decomposition is a spectral decomposition of the form $R = Q^T M Q$, where $R$ is diagonal and $Q$ is a real orthogonal matrix. By construction, the leading \texttt{eigenvalue\_count} diagonal entries of $R$ correspond to the desired eigenvalues, and the corresponding columns of $Q$ yield the associated eigenvectors.
\par
The requirement for symmetric input matrices stems from this structure of the \textsc{Schur} decomposition: if the input matrix is not symmetric, the resulting matrix $Q$ becomes a complex unitary matrix, and $R$ becomes tridiagonal. Extracting eigenvalues and eigenvectors from such a decomposition is considerably more involved. In fact, \textsc{ArnoldiMethod.jl} delegates this process to LAPACK routines, thereby restricting the computation to standard IEEE formats such as \texttt{float32} and \texttt{float64}, which is not feasible in this context.
\par
To obtain a high-precision reference solution for a given input matrix, we first apply the implicitly restarted \textsc{Arnoldi} method using \texttt{float128} arithmetic. This yields a set of \texttt{eigenvalue\_count} reference eigenvalues and a corresponding matrix whose columns contain the associated eigenvectors.
\par
When evaluating alternative numeric formats, directly applying the \textsc{Arnoldi} method to each target arithmetic in a naïve manner proves insufficient. This is due to the sensitivity of \textsc{Krylov} subspace methods to numerical perturbations, particularly when eigenvalues are tightly clustered. Small differences in arithmetic precision can cause these eigenvalues to appear in differing orders. Consequently, when comparing the computed eigenvectors to the reference solution, results may appear incorrect due solely to permutations, despite being otherwise accurate.
\par
We address this issue using a novel method, for which the authors have not identified any references in the existing literature. Instead of determining only the \texttt{eigenvalue\_count}, we compute an additional \texttt{eigenvalue\_buffer\_count} eigenvalues and eigenvectors, both during the initial preparation and throughout the evaluation. In all experiments we set this parameter to two. This strategy ensures that we avoid situations where a cluster of closely spaced eigenvalues might be cut off at the end of the desired eigenvector count, potentially leading to the omission of critical eigenvectors. By determining a few extra eigenvalues and eigenvectors, we provide additional \enquote{headroom} that can later be discarded, after the optimal permutation has been identified to match similar eigenvectors.
\par
To find this optimal permutation, we utilize the absolute cosine similarity matrix as a measure of similarity. Let $R := (\bm{r}_1 \cdots \bm{r}_m)$ represent the reference eigenvector matrix and $S := (\bm{s}_1 \cdots \bm{s}_m)$ the computed eigenvector matrix. The absolute cosine similarity matrix is defined as
\begin{equation}
    C_{i,j} := \frac{|\langle \bm{r}_i, \bm{s}_j \rangle|}{{\|\bm{r}_i\|}_2 {\|\bm{s}_j\|}_2}.
\end{equation}
This similarity measure yields positive values even when two vectors are opposites, as eigenvectors are only unique up to a sign.
\par
Having obtained the cosine similarity matrix, the next step is to find a permutation that yields the highest overall cosine similarity between each computed eigenvector and the corresponding reference eigenvector. To achieve this, we employ the Hungarian algorithm, which is conveniently available in Julia through the \texttt{Hungarian.jl} package. This algorithm takes the negative absolute cosine similarity matrix as input and returns the ideal permutation. Given the relatively small size of the absolute cosine similarity matrix, the Hungarian algorithm runs efficiently, despite its $\mathcal{O}(n^3)$ complexity. Once the permutation is determined, it is applied to the computed eigenvector matrix.
\par
As previously mentioned, eigenvectors are only unique up to a sign, meaning that the computed eigenvectors may be the exact opposite of the reference eigenvectors. Simply using the first entry of each vector as a sign reference is not stable, as the first entry may be very small, leading to potential sign flips at low precision. A more stable approach is to identify the index of the largest absolute entry in each reference eigenvector during preparation and use the sign of this entry in the computed eigenvector as the sign reference. This method ensures that the final computed eigenvector matrix is correctly signed, and the same permutation is applied to the eigenvalue vector.
\par
To quantify the errors, we only take the first \texttt{eigenvalue\_count} eigenvalues and the corresponding columns in the eigenvector matrix, both in the reference and the computed matrices. The absolute error is then computed by taking the $L^2$ norm of the difference between the reference and computed eigenvectors. The relative error is obtained by dividing the $L^2$ norm of the difference by the $L^2$ norm of the reference eigenvectors.
\section{Results}\label{sec:results}
In this section, we present the results of the experiments outlined previously. For each experiment, the computed relative errors are sorted and plotted as cumulative error distributions. This allows for a direct visualization of the proportion of runs that remain below a given relative error threshold. Moreover, because the datasets are diverse, sorting the errors still enables meaningful comparisons across different configurations.
\subsection{General Matrices}
\begin{figure*}[tbp]
	\begin{center}
        \subfloat[8 bits]{
			\includegraphics{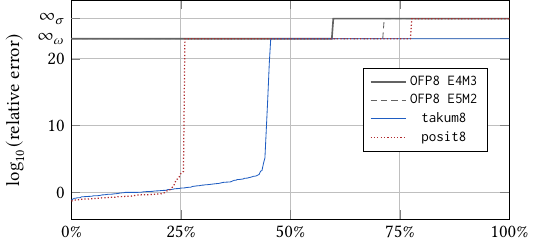}
			\includegraphics{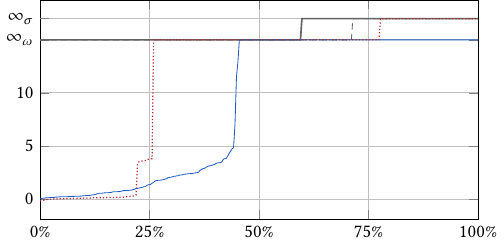}
        }\\
        \subfloat[16 bits]{
			\includegraphics{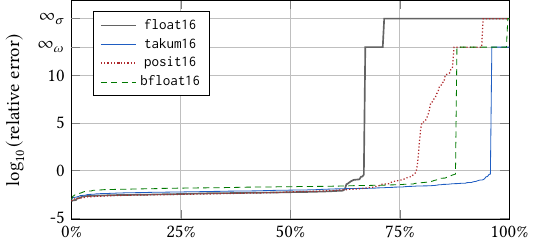}
			\includegraphics{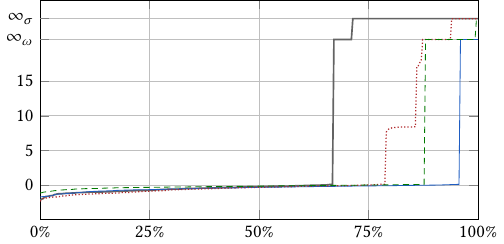}
        }\\
        \subfloat[32 bits]{
			\includegraphics{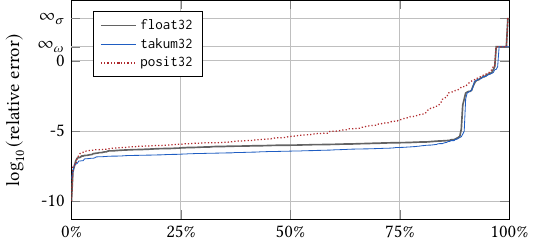}
			\includegraphics{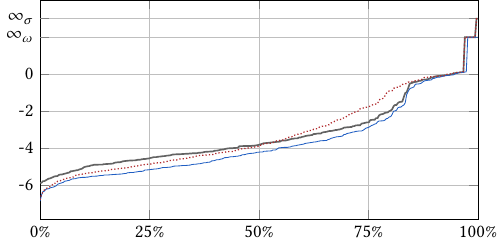}
        }\\
        \subfloat[64 bits]{
			\includegraphics{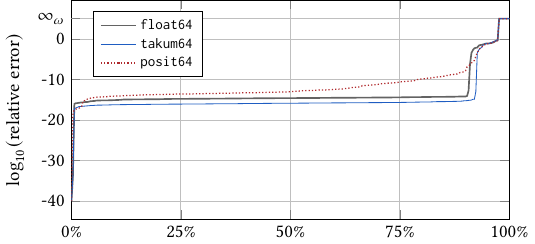}
			\includegraphics{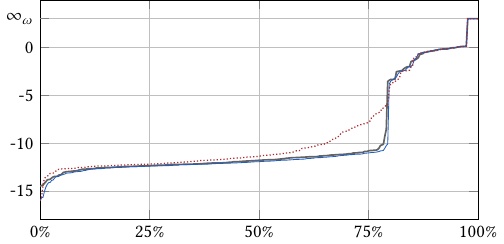}
        }
	\end{center}
	\caption{
        Cumulative error distribution of the relative errors of
        the $10$ largest eigenvalues (left) and their corresponding eigenvectors (right)
        of the general matrices computed using a range of machine number types.
        The symbol $\infty_\omega$ denotes where the \textsc{Arnoldi}
        method did not converge, $\infty_\sigma$ denotes where the dynamic
        range of the matrix entries exceeded the target number type.
	}
	\label{fig:general}
\end{figure*}
The results for the general matrices sourced from the SuiteSparse Matrix Collection are shown in Figure~\ref{fig:general}. Each row corresponds to experiments conducted with 8-, 16-, 32-, and 64-bit number formats, respectively. The left column illustrates the relative errors of the eigenvalues, while the right column displays the relative errors of the eigenvectors.
\par
At 8 bits, the OFP8 types consistently fail---either due to insufficient dynamic range or because the \textsc{Arnoldi} method fails to converge. While posits show slightly lower relative error in the lower percentiles, takums demonstrate greater overall stability.
\par
At 16 bits, both takums and \texttt{float16} achieve the lowest relative errors in the initial portion of the distribution, with takums performing slightly better overall. In contrast, \texttt{bfloat16} exhibits significantly higher relative errors but proves more stable than both \texttt{float16} and \texttt{posit16}. Notably, \texttt{takum16} emerges as the most stable format, outperforming even \texttt{bfloat16} in this regard.
\par
The 32-bit results reveal an unexpected trend: posit formats consistently underperform compared to \texttt{float32}, while takums achieve relative errors approximately half an order of magnitude lower than \texttt{float32}.
\par
This pattern becomes even more pronounced at 64 bits. Here, posits perform up to five orders of magnitude worse than \texttt{float64}. In contrast, takums consistently surpass \texttt{float64} by at least one order of magnitude in terms of eigenvalue relative error. For eigenvectors, takums perform comparably to \texttt{float32}, while posits again lag behind by up to five orders of magnitude.
\subsection{Biological Graphs}
\begin{figure*}[tbp]
	\begin{center}
        \subfloat[8 bits]{
			\includegraphics{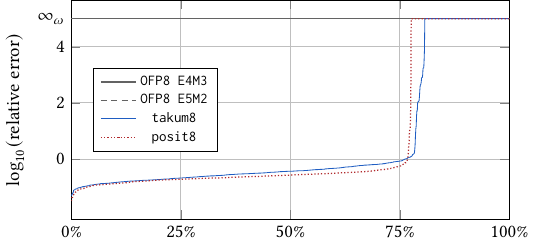}
			\includegraphics{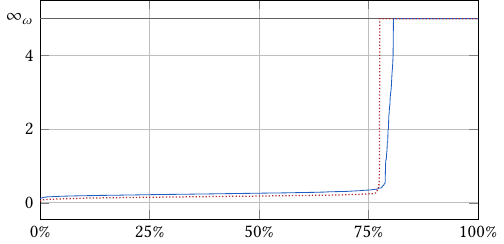}
        }\\
        \subfloat[16 bits]{
			\includegraphics{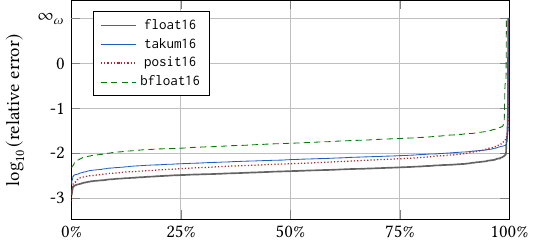}
			\includegraphics{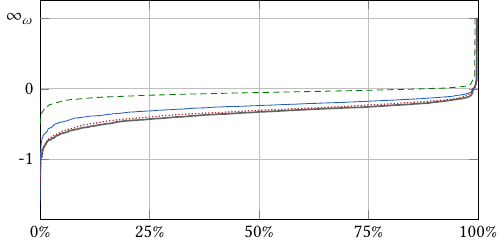}
        }\\
        \subfloat[32 bits]{
			\includegraphics{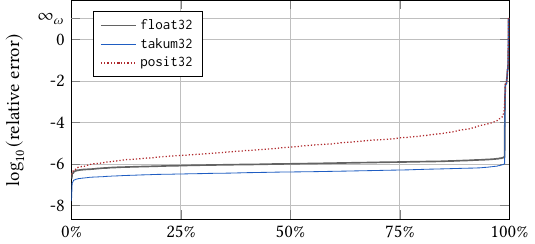}
			\includegraphics{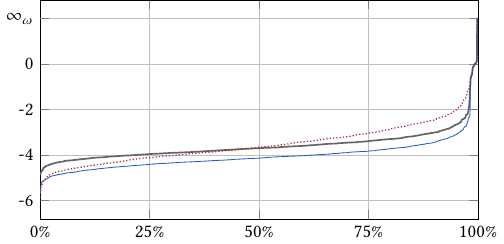}
        }\\
        \subfloat[64 bits]{
			\includegraphics{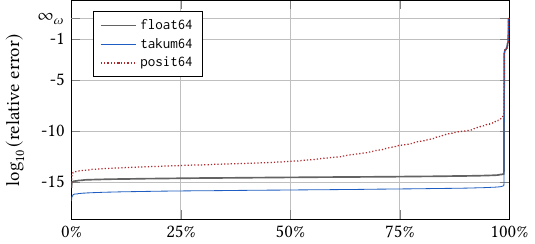}
			\includegraphics{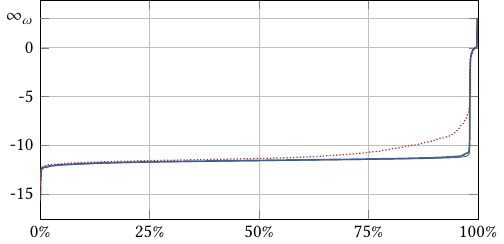}
        }
	\end{center}
	\caption{
        Cumulative error distribution of the relative errors of
        the $10$ largest eigenvalues (left) and their corresponding eigenvectors (right)
        of the biological graph symmetrized, normalized Laplacian
        matrices computed using a range of machine number types.
        The symbol $\infty_\omega$ denotes where the \textsc{Arnoldi}
        method did not converge.
	}
	\label{fig:biological}
\end{figure*}
We now turn to the graph matrices, beginning with the biological graphs, whose results are depicted in Figure~\ref{fig:biological}. At 8 and 16 bits, the results closely mirror those observed for the general matrices, with one notable exception: \texttt{bfloat16} exhibits a significantly higher relative error in this domain. Both posits and takums perform comparably to each other---worse than \texttt{float16}, but better than \texttt{bfloat16}. At 32 bits, posits are up to two orders of magnitude less accurate than \texttt{float32} in terms of eigenvalue errors, while takums consistently outperform \texttt{float32} by approximately half an order of magnitude across both eigenvalue and eigenvector errors. This trend becomes even more pronounced at 64 bits: takums surpass \texttt{float64} in both eigenvalue and eigenvector accuracy, whereas posits exhibit up to six and five orders of magnitude greater error in eigenvalues and eigenvectors, respectively.
\subsection{Infrastructure Graphs}
\begin{figure*}[tbp]
	\begin{center}
        \subfloat[8 bits]{
			\includegraphics{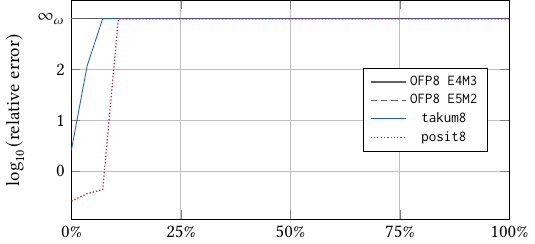}
			\includegraphics{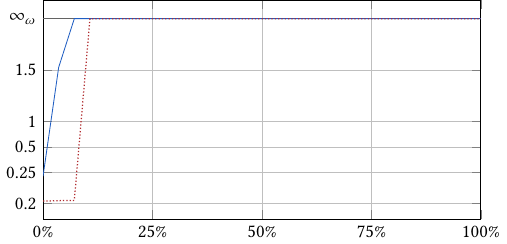}
        }\\
        \subfloat[16 bits]{
			\includegraphics{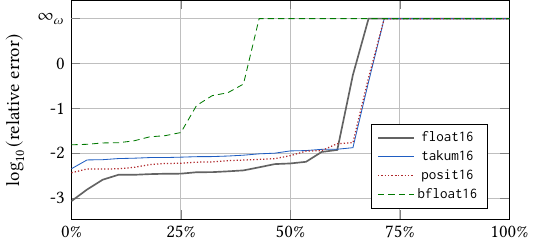}
			\includegraphics{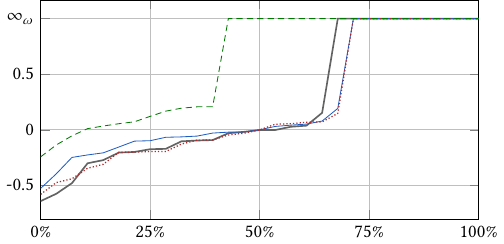}
        }\\
        \subfloat[32 bits]{
			\includegraphics{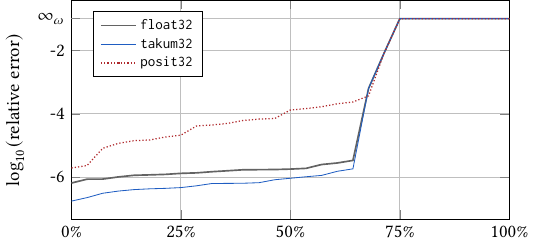}
			\includegraphics{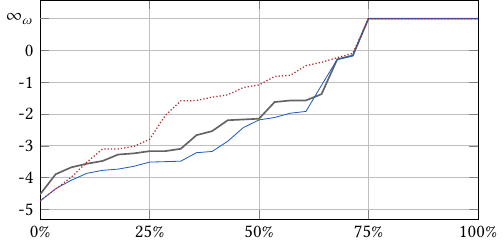}
        }\\
        \subfloat[64 bits]{
			\includegraphics{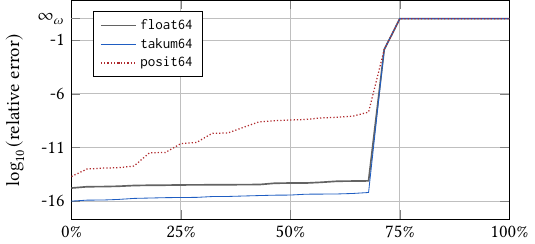}
			\includegraphics{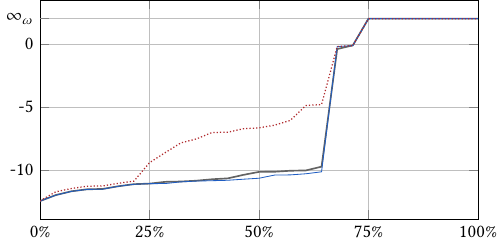}
        }
	\end{center}
	\caption{
        Cumulative error distribution of the relative errors of
        the $10$ largest eigenvalues (left) and their corresponding eigenvectors (right)
        of the infrastructure graph symmetrized, normalized Laplacian
        matrices computed using a range of machine number types.
        The symbol $\infty_\omega$ denotes where the \textsc{Arnoldi}
        method did not converge.
	}
	\label{fig:infrastructure}
\end{figure*}
The results for infrastructure graphs, shown in Figure~\ref{fig:infrastructure}, largely follow the same trends observed for biological graphs. At 8 bits, posits outperform takums. At 16 bits, both takums and posits lie between \texttt{float16} and \texttt{bfloat16}, although takums now exhibit slightly higher error than posits. At 32 and 64 bits, posits again display significantly higher errors than the corresponding IEEE 754 formats, whereas takums consistently outperform their IEEE counterparts in both eigenvalue and eigenvector accuracy.
\subsection{Social Graphs}
\begin{figure*}[tbp]
	\begin{center}
        \subfloat[8 bits]{
			\includegraphics{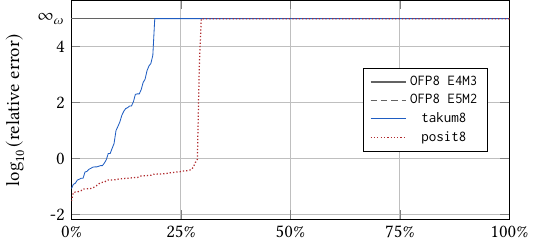}
			\includegraphics{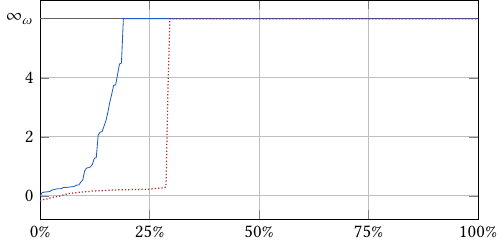}
        }\\
        \subfloat[16 bits]{
			\includegraphics{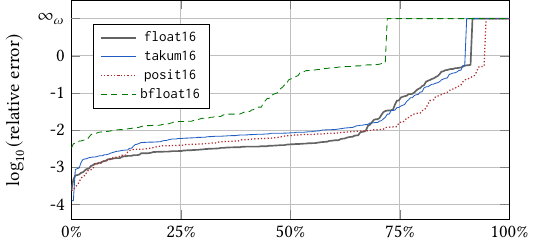}
			\includegraphics{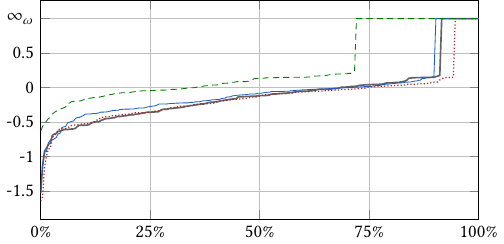}
        }\\
        \subfloat[32 bits]{
			\includegraphics{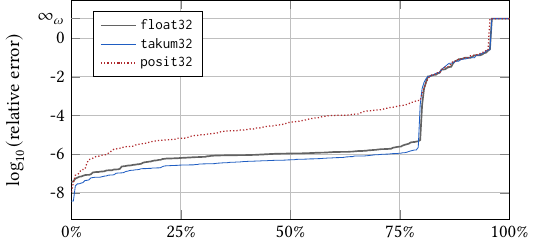}
			\includegraphics{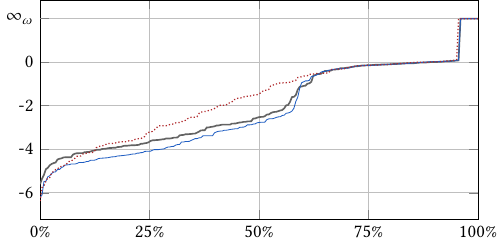}
        }\\
        \subfloat[64 bits]{
			\includegraphics{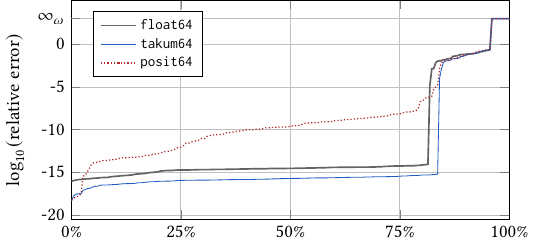}
			\includegraphics{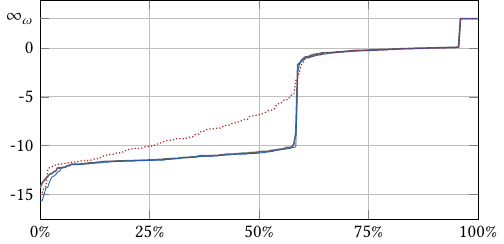}
        }
	\end{center}
	\caption{
        Cumulative error distribution of the relative errors of
        the $10$ largest eigenvalues (left) and their corresponding eigenvectors (right)
        of the social graph symmetrized, normalized Laplacian
        matrices computed using a range of machine number types.
        The symbol $\infty_\omega$ denotes where the \textsc{Arnoldi}
        method did not converge.
	}
	\label{fig:social}
\end{figure*}
The results for social graphs, presented in Figure~\ref{fig:social}, reinforce the overall performance trends observed in the previous graph categories. Notably, at lower precisions, the eigenvector errors of posits, takums, and \texttt{float16} are closely aligned. However, at 32 and 64 bits, posit formats demonstrate substantial degradation, with up to ten orders of magnitude greater error in eigenvalues and up to five orders of magnitude in eigenvectors, compared to the best-performing formats.
\subsection{Miscellaneous Graphs}
\begin{figure*}[tbp]
	\begin{center}
        \subfloat[8 bits]{
			\includegraphics{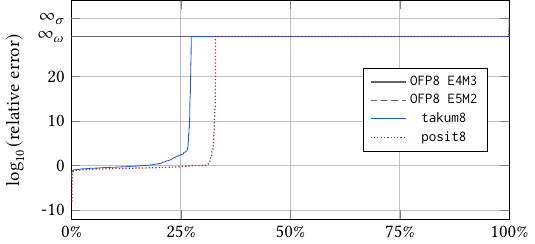}
			\includegraphics{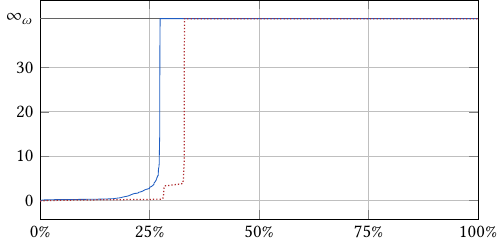}
        }\\
        \subfloat[16 bits]{
			\includegraphics{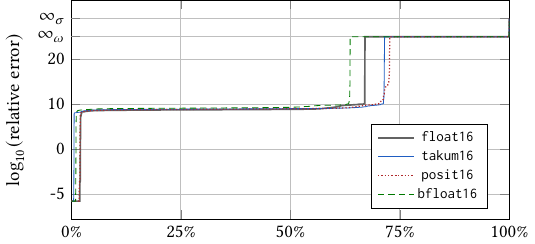}
			\includegraphics{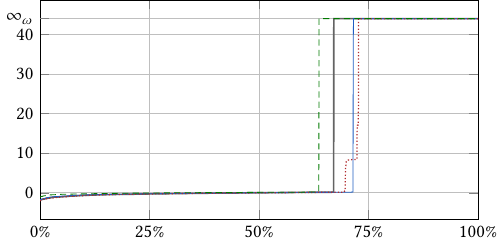}
        }\\
        \subfloat[32 bits]{
			\includegraphics{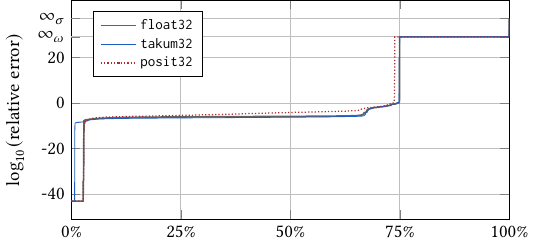}
			\includegraphics{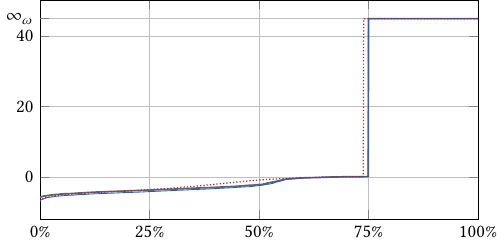}
        }\\
        \subfloat[64 bits]{
			\includegraphics{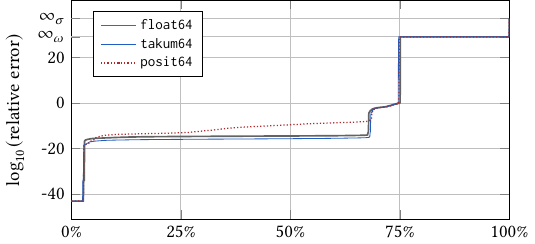}
			\includegraphics{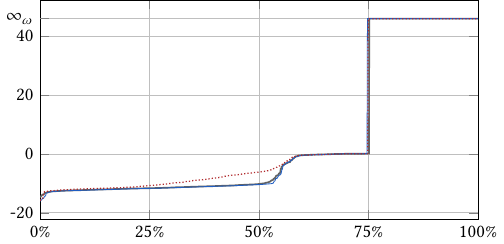}
        }
	\end{center}
	\caption{
        Cumulative error distribution of the relative errors of
        the $10$ largest eigenvalues (left) and their corresponding eigenvectors (right)
        of the miscellaneous graph symmetrized, normalized Laplacian
        matrices computed using a range of machine number types.
        The symbol $\infty_\omega$ denotes where the \textsc{Arnoldi}
        method did not converge, $\infty_\sigma$ denotes where the dynamic
        range of the matrix entries exceeded the target number type.
	}
	\label{fig:misc}
\end{figure*}
The results for the miscellaneous graphs are shown in Figure~\ref{fig:misc}. A distinguishing observation here is that none of the 16-bit formats manage to return useful results, suggesting that these matrices pose greater analytical challenges. At 32 and 64 bits, although the excessive error associated with posits remains present, it is somewhat less pronounced than in other graph categories. Nonetheless, takums continue to outperform both \texttt{float32} and \texttt{float64} across all metrics.
\section{Conclusion}\label{sec:conclusions}
This paper evaluated the number formats OFP8 E4M3 and E5M2, \texttt{bfloat16}, posits, and takums against IEEE 754 floating-point numbers within the \textsc{Arnoldi} method for computing select eigenvalues. The datasets comprised both general matrices and diverse graph Laplacians. It was demonstrated that the OFP8 formats are, overall, unsuitable for this purpose, and that \texttt{bfloat16} performs worse than \texttt{float16} in this application. Both tapered-precision formats evaluated in this study---\texttt{posit16} and \texttt{takum16}---yielded performance between that of \texttt{bfloat16} and \texttt{float16}.
\par
A particularly notable finding is the performance of posits at 32 and 64 bits, which resulted in relative errors in eigenvalues and eigenvectors several orders of magnitude higher than the corresponding IEEE 754 types. In contrast, takums consistently outperformed their IEEE counterparts. To the best of the author’s knowledge, this is one of the first documented cases where posits underperform to such a significant extent, marking these as important and novel results.
\par
Overall, the results indicate that takums, due to their high dynamic range and stable behavior across matrix types, are better suited for general-purpose numerical computations than posits.
\par
Future work could expand the evaluation to larger and more varied datasets, as well as explore downstream graph-based applications that rely on accurate eigenvalue and eigenvector computation.
\section*{Author Contributions}
Specific roles and contributions will be detailed following the review process.
\printbibliography

\end{document}